\documentclass[a4paper,10pt]{article}
\usepackage{stmaryrd}
\usepackage{bbm}
\usepackage{amsfonts}
\usepackage{amsmath}
\usepackage{amssymb}
\usepackage{mathrsfs}
\usepackage{accents}
\usepackage{amscd}
\usepackage{tikz}

\newtheorem{atheorem}{\bf \temp}[section]
\newtheorem{thm}[atheorem]{Theorem}
\newtheorem{cor}[atheorem]{Corollary}
\newtheorem{lem}[atheorem]{Lemma}
\newtheorem{prop}[atheorem]{Proposition}

\newtheorem{rem}[atheorem]{Remark}

\numberwithin{equation}{section}

\title{\textbf{On the  Steady Magnetohydrodynamic Equations with Nonhomogeneous Boundary Conditions }}
\author{ Xixia Ma  \footnote{Corresponding author.. E-mail addresses: kfmaxixia@163.com(Xixia Ma)} \ \ \ \ \ \ \
 \\}
\date{}

\begin{document}

\maketitle

\textbf{Abstract.}We study both  the topological structure stability  and the relations  of  the steady Magnetohydrodynamic equations  when $\nu,\eta$ are given different values  in muti-connected  bounded domain.  We also show the solutions's existence  for fixed $\nu, \eta.$  The theoretical  is the Morse-Sard theorem on Sobolev spaces.
\begin{center}
\item\section{Introduction}
\end{center}

  The Magnetohydrodynamics (MHD) is the macroscopic theory  of electrically conducting fluids, providing a powerful and practical theoretical framework for describing both laboratory and astrophysical plasmas, and it also describes the macroscopic behavior of electrically conducting fluids, notably of plasmas. However, here we  focus on  the interaction between magnetic fields and moving, conducting fluids. The mutual interaction of a magnetic field $\mathbf{B}$ and a velocity $\mathbf{u}$  arises partially as a result of the laws of Faraday and Amp$\grave{e}$re, and partially because of the Lorentz force experienced by a current-carrying body.

  In contrast to do with dynamics in MHD in the past,  we try to deal with magnetostatic configurations. It is based on the powerful mathematical theory on Sard theorem on Sobolev spaces. And in the past the static magnetic field configuration may appear in a bewildering variety of shapes generated by the fluid self, not by the particular boundary. And in this paper we focus on the nonlinear stability theory from the topological aspects and more consideration on this nonlinear phenomena from both the physics and the mathematics.

    Consider  the stationary Magnetohydrodynamic system in a bounded domain with type of (1.1), i.e. $\Omega$ denotes a bounded region with  $C^{2}$-smooth boundary $\partial\Omega=\cup^{N}_{j=0}\Gamma_{j}$ consisting of N+1 disjoint components $\Gamma_{j};$
\begin {equation}
\Omega=\Omega_{0}\setminus(\cup^{N}_{j=0}\bar{\Omega}_{j}),\bar{\Omega}_{j}\subset\Omega_{0},j=1,\ldots,N,
\end {equation}
where $\Gamma_{j}=\partial\Omega_{j}$,\footnote{Here we assume this condition is the plasma boundary condition,that is, $\mathbf{B}\cdot \mathbf{n}=0,\mathbf{B}\times \mathbf{n}=\mathbf{b}(x),$ where $\mathbf{n}$ is the outer normal vector of $\partial\Omega.$ Because $\mathbf{B}$ is a 2-dim vector field,$(1.1)_{6}$ is well-defined,see also [2].In this paper, without generality of loss, we assume that $\mathbf{a}\cdot\mathbf{b}=0.$}
and
\begin{equation}
\left\{\begin {array}{l}
-\nu\triangle\mathbf{ u}+\mathbf{u}\cdot\nabla \mathbf{u}+\nabla p=\mathbf{B}\cdot\nabla \mathbf{B},\\
-\eta\triangle \mathbf{B}+\mathbf{u}\cdot\nabla \mathbf{B}=\mathbf{B}\cdot\nabla \mathbf{u},\\
\nabla\cdot \mathbf{u}=0, \nabla\cdot \mathbf{B}=0,\\

\mathbf{ u}(x)|_{\partial\Omega}=\mathbf{a}(x), \mathbf{B}(x)|_{\partial\Omega}=\mathbf{b}(x),
\end{array}\right.
\end{equation}
where $\mathbf{B}$ is a magnetic field, and $\mathbf{u}$ is a velocity, $\eta,\nu,$ the magnitude dissipation coefficient and viscous coefficient, respectively.
The continuity equation $(1.2)_{3} $ implies  the necessary compatibility condition for the solvability of the problem (1.2):
\begin {equation}
\int_{\partial\Omega} \mathbf{a}\cdot\mathbf{ n}dS=\sum^{N}_{j=0}\mathcal{F}_{j}=\int_{\partial\Omega} \mathbf{b}\cdot \mathbf{n}dS=\sum^{N}_{j=0}\mathcal{G}_{j}=0,
\end {equation}
where $\mathbf{n}$ is a unit vector of the outward ( with respect to $\Omega$) normal to $\partial\Omega$ and $\mathcal{F}_{j}=\int_{\Gamma_{j}}\mathbf{a}\cdot\mathbf{ n} dS,\quad \mathcal{G}_{j}=\int_{\Gamma_{j}}\mathbf{b}\cdot \mathbf{n} dS.$

The equality (1.3) implies that the total flux through the boundary is zero and  the total magnetic flux through the boundary is zero. And we assume $\nu=\eta,$  which means dissipation contributions both the magnetics and the fluid.

If one replaces (1.3) by the stronger condition
\begin {equation}
\int_{\Gamma_{j}}\mathbf{a}\cdot \mathbf{n} dS=\int_{\Gamma_{j}}\mathbf{b}\cdot \mathbf{n} dS=0,j=0,\ldots,N,
\end {equation}
The existence of a weak solution $(\mathbf{u},\mathbf{B})\in W^{1,2}(\Omega)$ to the problem(1.1),(1.2) is trivial only under assumption(1.4), or even for sufficiently small fluxes $\mathcal{F}^{2}_{j},\mathcal{G}^{2}_{j},$ or under certain conditions on the domain $\Omega$ and the boundary value $\mathbf{a},\mathbf{b}$ in [11].

We only assume throughout this paper that  the equalities (1.3) hold.  Our intention is to show some topological property and the relations  of the MHD system (1.2) under different values of $\nu,\eta$ in the above domain $\Omega$ for some $\mathbf{a},\mathbf{b}$ belonging to some space and the existence of the solution of the  MHD system for fixed $\nu,\eta.$

We change the system (1.2) in the following transformation,  let
$$\mathbf{Z}^{+}=\mathbf{u}+\mathbf{B},\quad \quad \mathbf{Z}^{-}=\mathbf{u}-\mathbf{B},$$
then (1.2) becomes
\begin{equation}
\left\{\begin {array}{l}
-\nu\triangle \mathbf{Z}^{+}+\mathbf{Z}^{-}\cdot\nabla \mathbf{Z}^{+}+\nabla p=0,\\
-\nu\triangle \mathbf{Z}^{-}+\mathbf{Z}^{+}\cdot\nabla \mathbf{Z}^{-}+\nabla p=0,\\
\nabla\cdot \mathbf{Z}^{+} =0, \nabla\cdot \mathbf{Z}^{-}=0,\\

 \mathbf{Z}^{+}(x)|_{\partial\Omega}=\mathbf{a}(x)+\mathbf{b}(x), \mathbf{Z}^{-}(x)|_{\partial\Omega}=\mathbf{a}(x)-\mathbf{b}(x).
\end{array}\right.
\end{equation}

The main result of this paper is the following theorem to show of the solution of (1.1)-(1.4).

\begin {thm} Let $\Omega$ be a bounded domain with  $C^{2}$-smooth boundary $\partial\Omega=\cup^{N}_{j=0}\Gamma_{j}$ consisting of N+1 disjoint components $\Gamma_{j};$ i.e.,
$$
\Omega=\Omega_{0}\setminus(\cup^{N}_{j=0}\bar{\Omega}_{j}),\bar{\Omega}_{j}\subset\Omega_{0},j=1,\ldots,N,
$$
where $\Gamma_{j}=\partial\Omega_{j}$ .

If $\mathbf{a},\mathbf{b}\in W^{\frac{3}{2},2}(\partial\Omega)$ and let conditions(1.3) be fulfilled. When $\nu=\eta=0,$ then  the ideal MHD satisfies some Bernoulli's law.
\end {thm}
\begin{thm}
Under the conditions of Theorem 1.1,  when $\nu=\eta\rightarrow0,$  the MHD system preserve the stability on the Bernoulli's law of Theorem 1.1.
\end{thm}
\begin{thm}
Under the above conditions, then the problem (1.2) admits at least one weak solution $(\mathbf{u},\mathbf{B})$ in $W^{1,2}(\Omega)\times W^{1,2}(\Omega),$ for fixed $\nu=\eta>0.$
\end{thm}
\begin {rem}
It is well known that under the hypotheses of Theorem 1.3, every weak solution $(\mathbf{u},\mathbf{B})$ of the problem(1.2) is more regular:$(\mathbf{u},\mathbf{B})\in(W^{2,2}(\Omega)$
$\cap W_{loc}^{3,2}(\Omega),W^{2,2}(\Omega)$
$\cap W_{loc}^{3,2}(\Omega)).$ Generally speaking, the solution is as regular as the data allow.
\end {rem}

Another classical model for the viscous incompressible steady fluid is the Navier-Stokes model. the nonhomogeneous boudary value problem for Navier-Stokes equations in any bounded domain in $\mathbb{R}^{2} $ with multiply connected boundary was studied by several authors in [6,9] and so on.

This paper is organized as follows. Section 2 is a preliminary and Morse-Sard theorem.  In Sections 3-5, we give the corresponding proofs of Theorems 1-3,respectively.
\begin{center}
\item\section{Preliminary and Morse-Sard theorem on Sobolev spaces}
\end{center}
$\mathbf{2.1.}$ Some simple properties  of the low-dimensional  topology.

\begin {lem}(see[10])  Let $\Omega\subset\mathbb{R}^{2}$ be a bounded domain with Lipschitz boundary and let $K\subset\bar{\Omega}$ be a continuum. Then there exists $\delta>0$ such that for any continuous injective function $\gamma:I=[0,1]\rightarrow\bar{\Omega}$ with the properties $\gamma(0),\gamma(1)\in K,$ and $\gamma((0,1))\subset\{x\in\Omega: dist(x,K)<\delta\},$ the following assertion is valid.

 $(\mathfrak{N})$ For any interval $(\alpha,\beta)$ adjoining the set $\tilde{I}=\gamma^{-1}(K)$ (that is, $\tilde{I}$ is a compact subset of the interval [0,1] and $0,1\in\tilde{I},(\alpha,\beta)$is a connected component subset of the open set $(0,1)\setminus\tilde{I})$ there exists a continuum $K_{\alpha\beta}\subset K$ and a simply connected domain $\Omega_{\alpha\beta}\subset\Omega$ such that $\Omega_{\alpha\beta}\cap K=\emptyset,\gamma(\alpha),\gamma(\beta)\in K_{\alpha\beta}$ and $\partial\Omega_{\alpha\beta}=K_{\alpha\beta}\cup\gamma([\alpha,\beta]).$
 \end {lem}
$\mathbf{ 2.2.}$ The properties on Sobolev functions and Morse-Sard theorem on Sobolev spaces

 By a domain we mean an open connected set.In this paper we deal with bounded domain $\Omega\subset\mathbb{R}^{2}$ be a bounded domain with type(1.1)

 Here we follow the notations in [13]. We still denote $ H(\Omega)$ is a subspace of all solenoidal vector fields $(\nabla\cdot u=0)$ from $\dot{H}^{1}(\Omega)$ with the norm $\|u\|_{H(\Omega)}=\|\nabla u\|_{L^{2}(\Omega)}.$

 In the following, we state several results about Sobolev functions

\begin {lem} (see[6]) If $w\in W^{1,s}(\mathbb{R}^{2}),s\geq1,$ then exists a set $A_{1,w}\subset\mathbb{R}^{2}$ with the following properties:

(i) $ \mathfrak{h}^{1}(A_{1,w})=0;$

(ii) for each $x\in\Omega\setminus A_{1,w},$
$$\lim_{r\rightarrow0}\frac{1}{|B_{r}(x)|}\int_{B_{r}(x)}|w(z)-w(x)|^{2}dz=0;$$

(iii) for every $\varepsilon>0,$ there exists a set $U\subset\mathbb{R}^{2}$ with $ \mathfrak{h}^{1}_{\infty}(U)<\varepsilon $ and $A_{1,w}\subset U $ such that the function $ w $ is continuous on $\bar{\Omega}\setminus U;$

(iv) for every unit vector $\mathbf{I}\in\partial B_{1}(0)$ and almost all straight lines $L$ parallel to $\mathbf{I},$ the restrictions $w|_{L}$ is an absolutely continuous function (of one variable).
\end {lem}

Here and henceforth we denote by $\mathfrak{h}^{1}$ the one-dimensional Hausdorff measure, i.e., $ \mathfrak{h}^{1}(F)=\lim_{t\rightarrow0^{+}}\mathfrak{h}^{1}_{t}(F),$ where
$$\mathfrak{h}^{1}_{t}(F)=\inf\bigg\{{\sum^{\infty}_{i=1}diamF_{i}: diamF_{i}\leq t, F\subset\cup^{\infty}_{i=1}F_{i}}\bigg\}.$$

 The theorem below is due to J.Bourgain,M.Korobkov and J.Kristensen[3]. It is the Morse-Sard theorem on Sobolev spaces.

 \begin {thm}
 Let $\Omega\subset\mathbb{R}^{2}$ be a bounded domain with Lipschitz boundary. If $\psi\in W^{2,1}(\Omega),$ then

 (i) $\mathfrak{h}^{1}(\{\psi(x):x\in\bar{\Omega}\setminus A_{\psi},\nabla\psi(x)=0\}).$

 (ii)  For every $\varepsilon>0,$ there exists $\delta>0$ such that $\mathfrak{h}^{1}(\psi(U))<\varepsilon$ for any set $U\subset\bar{\Omega}$ with $ \mathfrak{h}^{1}_{\infty}(U)<\delta.$

 (iii) For every $\varepsilon>0,$ there exist an open set $V\subset\mathbb{R}$ with $\mathfrak{h}^{1}(V)<\varepsilon$ and a function $g\in C^{1}(\mathbb{R}^{2})$ such that for each $x\in\bar{\Omega}$ if $\psi(x)\notin V,$ then $x\notin A_{\psi}$
 and $\psi(x)=g(x),\nabla\psi(x)=\nabla g(x)\neq0.$

 (iv)  For $\mathfrak{h}^{1}-$ almost all $y\in\psi(\Omega)\subset\mathbb{R},$ the preimage $\psi^{-1}(y)$ is a finite disjoint family of $C^{1}-$ curves $S_{j},j=1,2,\ldots,N(y).$ Each $ S_{j}$ is either a cycle in $\Omega$ (i.e.,$ S_{j}\subset\Omega $is homeomorphic to the unit circle $ \mathbb{S}^{1}$) or a simple arc with endpoints on $\partial\Omega.$(In this case $S_{j}$is transversal to $\partial\Omega.$)
 \end {thm}

\begin{rem}  In my opinion, the above theorem is not only the generalization of the classical Morse-Sard theorem, but also the bridge between the classical PDE on Sobolev spaces and the modern geometry from the intrinsic property of the geometry. And this theorem is the theoretical basis in  my paper.
\end{rem}

\begin {rem}(see[12])  If $f\in W^{2,1}(Q),$ then by Theorem 2.3(iv), there exists a dense subset $E$ of $(0,1)$ such that $\varphi(t)$ is a $C^{1}-$ curve for every $t\in E.$ Moreover, $\varphi(t)$ is either a cycle or a simple arc with endpoints on $\partial Q.$
\end {rem}
$ \mathbf{2.3.}$ Some properties on continuous  functions

\begin {lem}
Let $f\in C(Q),$ then for any two different points $A,B\in T_{f},$ there exists an injective function $\varphi:[0,1]\rightarrow T_{f}$ such that

(i) $\varphi(0)=A, \varphi(1)=B;$

(ii) for any $t_{0}\in[0,1],$ the convergence $lim_{[0,1]\ni t\rightarrow t_{0}}sup_{x\in\varphi(t)}dist(x,\varphi(t_{0}))\rightarrow0$ holds.

(iii) for any $t_{0}\in(0,1),$ the sets $ A,B$ lie in the different connected components of the set $ Q\setminus \varphi(t).$

\end {lem}

\begin {rem}(see[12])
All results of Lemmas 2.6 remain valid for level sets of continuous functions $f:\bar{\Omega}\rightarrow\mathbb{R},$ where $\Omega$ is a multi-connected bounded domain of type(1.1), provided $f\equiv\xi_{j}= const$ on each inner boundary component $\Gamma_{j},j=1,\ldots,N.$ Indeed, we can extend $ f $ to the whole $\bar{\Omega}_{0}$ by putting $f\equiv\xi_{j}$
for $x\in\bar{\Omega}_{j}, j=1,\ldots,N.$  The extended function $f$ will be continuous on the set $\bar{\Omega}_{0}$ that is homeomorphic to the unit square $Q=[0,1]^{2}.$
\end {rem}

\begin{center}
\item\section{Some properties of the ideal MHD equations }
\end{center}

Now we consider the ideal MHD system in $\Omega$ in the following form,

\begin{equation}
\left\{\begin {array}{l}
\mathbf{v}\cdot\nabla\mathbf{ v} + \mathbf{B}\cdot\nabla \mathbf{B}+\nabla p=0,\\
\mathbf{v}\cdot\nabla \mathbf{B}-\mathbf{B}\cdot\nabla \mathbf{v}=0,\\
\nabla\cdot \mathbf{v}=\nabla\cdot \mathbf{B}=0,\\
\mathbf{v}|_{\partial\Omega}=\mathbf{a}(x),\mathbf{B}|_{\partial\Omega}=\mathbf{b}(x),
\end{array}\right.
\end{equation}
where $\mathbf{a}(x), \mathbf{b}(x)$ may be constants. Here we also that $\mathbf{a}\cdot\mathbf{b}=0.$

Using the transformation $\mathbf{w}_{1}=\mathbf{v}-\mathbf{B},\mathbf{w}_{2}=\mathbf{v}+\mathbf{B},$
\begin{equation}
\left\{\begin {array}{l}
(\mathbf{w}_{1}+2\mathbf{B})\cdot\nabla \mathbf{w}_{1} +\nabla p=0,\\
(\mathbf{w}_{2}-2\mathbf{B})\cdot\nabla \mathbf{w}_{2} +\nabla p=0,\\
\nabla\cdot\mathbf{ w}_{1}=\nabla\cdot\mathbf{ w}_{2}=\nabla\cdot \mathbf{B}=0,\\
\mathbf{w}_{1}|_{\partial\Omega}=\mathbf{a}-\mathbf{b},\mathbf{w}_{2}|_{\partial\Omega}=\mathbf{a}+\mathbf{b},\mathbf{B}|_{\partial\Omega}=\mathbf{b}.
\end{array}\right.
\end{equation}

Let $\Omega\subset\mathbb{R}^{2}$ be an bounded domain with Lipschitz boundary. Assume that $\mathbf{v},\mathbf{B}\in W^{1,2}(\Omega)$ and $p\in W^{1,s}(\Omega),s\in[1,2)$ satisfy the equations (3.1) for almost all $ x\in\Omega.$  Since $\nabla\cdot\mathbf{ v}=\nabla\cdot \mathbf{B}=0,$ we can introduce both stream function and magnetic stream function  $\psi_{1},\psi_{2}\in W^{2,2}(\Omega).$ Denote $\psi=\psi_{1}\pm\psi_{2},\nabla\psi_{1}=(-v_{2},v_{1}),\nabla\psi_{2}=(-B_{2},B_{1}).$  Denote by $\Phi=p+\frac{1}{2}|\mathbf{v}\pm\mathbf{B}|^{2},$  then $\Phi$ can be regarded as total energy per unit mass of the system (3.1). Obviously, $\Phi\in W^{1,s}(\Omega),$ for all $s\in[1,2).$ By direct calculations, one easily gets the identity
\begin {equation}
\nabla\Phi\equiv[\triangle\psi_{1}\pm\left(
\begin {array}{l l}
0 \quad 1\\
-1 \quad 0
\end {array}
\right)
\nabla^{2}\psi_{2}
]\nabla\psi.
\end {equation}

\begin{rem} In fact, according to Amp$\grave{e}$re's law, the induced currents must give rise to ta second, induced magnetic field, since the boundary values vanish. In addition, the induced magnetic field interacts with the induced current and generates a Lorentz force, which will inhibit the relative movement of the magnetic field and the fluid.
\end{rem}
\begin{rem} In this paper, we mainly concern the geometric structure of the fluid and the magnetic field. From the definition of $\psi_{1},\psi_{2}$  and  the equality (3.3), we can also see the mutual interaction of a magnetic field $\mathbf{B}$ and a velocity $\mathbf{v},$  which inhibit each other.  It is worth stating now, from the Euler potential representation of the magnetic field $\mathbf{B} $ in multiply connected ($\mathbf{B}=\nabla\times\mathbf{A},$ where $\mathbf{A}=\alpha\nabla\beta+\nabla\chi$), as for instance in a torus, the representation may not be single-valued, then the corresponding magnetic helicities is not gauge-invariant, which is the key difference between the simple-connected domain and the muilti-connected domain from the geometric viewpoint.
\end{rem}

Now we state the main result in this section. In this section,we prove the case$(\mathbf{v}-\mathbf{B})$. The result implies  the geometric structure distribution  of the total energy  in any subdonmain of $\Omega.$

\begin {thm}
Let $\Omega\subset\mathbb{R}^{2}$ be a bounded multiply connected domain with $C^{2}$ boundary $\partial\Omega=\cup^{N}_{j=1}\Gamma_{j}.$ Assume that  $\mathbf{v},\mathbf{B}\in W^{1,2}(\Omega)$ and $p\in W^{1,s}(\Omega),s\in[1,2)$ satisfy the the ideal MHD equations (3.1) for almost all $ x\in\Omega.$ Then for any connected set $ K\subset\bar{\Omega}$ such that
\begin {equation}
\psi|_{K}=const,
\end {equation}
the assertion
there exists a constant $C=C(K)$ such that
\begin {equation}
\Phi(x)=C
\end {equation}
for $\mathfrak{h}^{1}-$ almost all $x\in K$ holds.
\end {thm}
\textbf{Proof.}  we will decompose the proof into two steps.

Step 1, It is trivial that for almost all $y\in\psi(\Omega)$ and any connected component $K$ of $\psi^{-1}(y),$ the theorem holds. At the same time, the preimage of the corresponding $\Phi(x)=const$ is a disjoint finite $C^{1}$ circle $S(j),j=1,\ldots,N(y),$ according to Theorem 2.3(iv).

Step 2,  we need to prove for any $y\in\psi(\bar{\Omega}),$ and any connected component $K$ of $\psi^{-1}(y),$ the theorem is also valid.

In fact,  choose $A_{\psi}$ in Lemma 2.2 (i) and let $a,b \in K\setminus A_{\psi},$   we need to prove that $\Phi(a)=\Phi(b).$
 From Lemma 2.2 (iv), we can choose a function $\zeta(s)$ satisfies the following properties:
 \begin{itemize}
 \item[{i.}]  $\zeta(s):[0,1]\rightarrow \bar{\Omega}\setminus A_{\psi}$ is a Lipschitz and injective function;
 \item[{ii.}]  $\zeta(0)=a,\zeta(1)=b, \tilde{K}=\zeta^{-1}(K)\subset[0,1];$
 \item[{iii.}]  for an arbitrary interval $(\alpha,\beta),$ if $(\alpha,\beta)$ is a connected component of the open set $(0,1)\setminus\tilde{K},$  then there is a compact and simple connected set $K_{\alpha\beta}\subset K$ and a simply connected subdomain $\Omega_{\alpha\beta}\subset\Omega$ such that $\partial\Omega_{\alpha\beta}=K_{\alpha\beta}\cup \zeta([\alpha,\beta]).$
 \end{itemize}

 In the subdomain $\Omega_{\alpha\beta},$  we know the fact that the collection of connected components of level sets of $\psi|_{\bar{\Omega}_{\alpha\beta}}$ equipped with the natural topology is a topological space, which is homeomorphic to a tree $T.$

 From Lemma 2.6, $\Omega_{\alpha\beta}$ simple connected  and $\partial\Omega_{\alpha\beta}=K_{\alpha\beta}\cup \zeta([\alpha,\beta]),$  there exist the sequences $t_{i}\rightarrow\alpha+0$ and $s_{i}\rightarrow\beta-0$ such that $\zeta(t_{i})$ and $\zeta(s_{i})$ can be jointed by a smooth arc from  the level set of $\psi|_{\Omega_{\alpha\beta}}$ on which $\Phi=const.$ And this implies that for any interval $(m,n)\subset(0,1),$ if $(m,n)$ contains only a finite set of points of $\tilde{K},$ $\int^{n}_{m}\Phi'(\zeta(s))ds=0.$

 Now we consider the closed set $K_{inf}=\{t\in[0,1]: $in any neighborhood of t, there are infinitely many points of $\tilde{K}\}.$  It is easy to check that
 $\int_{[0,1]\setminus K_{inf}}\Phi'(\zeta(s))ds$
 $=0.$

 It is easy to observe that  Lemma 2.2 implies that $\zeta(s)$ is differentiable for almost all $s\in[0,1].$ Therefore, if for $t\in K_{inf},$ there exists $\zeta'(s)$ and $\zeta'(s)\cdot\nabla\psi(\zeta(s))=0.$  We immediately derive $\int_{K_{inf}}\Phi'(\zeta(s))ds=0,$ then $\Phi(b)-\Phi(a)=\int_{[0,1]}\Phi'(\zeta(s)ds=0, $   Theorem 3.3 is proved.

\begin {rem}
In particular, if $\mathbf{v}=\mathbf{B}=constant$ on $\partial\Omega$(in the sense of trace), then the pressure $p(x)$ is constants on  $\partial\Omega$. Note that $p(x)$ could take different constant values on different connected components of the boundary. Furthermore,
$$p\in C(\bar{\Omega})\cap W^{1,2}(\Omega)$$  holds.
\end {rem}
\begin {lem}
If the conditions of Theorem 3.3 are satisfied, then there exist constants $\xi_{0},\ldots,\xi_{N}\in\mathbb{R}$ such that
$\psi(x)=\xi_{j}$ on each component $\Gamma_{j},j=0,\ldots,N.$
\end {lem}
\begin{rem} If the boundary consists of a disjoint finite $C^{2}-$smooth cycle, since $\mathbf{a}\cdot\mathbf{b}=0,$ without generality of loss, we assume $\mathbf{a}=\mathbf{b}=constant,$ from Remark 3.4 and Lemma 3.5,  then $\Phi(x)=C_{\Gamma_{j}},j=1,\ldots,N.$ However, although nearby each  boundary connected component $\Gamma_{j}$, from Theorem 3.3, there is a circle $S$ such that $\psi_{1}-\psi_{2}=constant,\psi_{1}+\psi_{2}=constant,$ then $p+\frac{1}{2}|\mathbf{v}-\mathbf{B}|^{2}=constant,$
$p+\frac{1}{2}|\mathbf{v}+\mathbf{B}|^{2}=constant.$  In general, it is obtained that $\mathbf{v}\cdot\mathbf{B}=constant\neq0,$ but $\mathbf{v}\cdot\mathbf{B}=0$ on $\Gamma_{j}.$ This change is regarded as the instability in the stability case, because it almost don't change the topology of the energy distribution, which implies the geometric structure stability.
\end{rem}

For $x\in\bar{\Omega},$ denote by $ K_{x}$ the connected component of the level set $\{z\in\bar{\Omega}:\psi(z)=\psi(x)\}$
containing the point $ x.$ By Lemma 3.5, $ K_{x}\cap\partial\Omega=\emptyset$ for every $y\in\psi(\bar{\Omega})\setminus\{\xi_{0},\ldots,\xi_{N}\}$ and for every $x\in\psi^{-1}(y).$ Thus, Theorem 2.3(ii) and (iv)
imply that for almost all $y\in\psi(\bar{\Omega})$ and for every $x\in\psi^{-1}(y),$ the equality $K_{x}\cap A_{\psi}=\emptyset$
holds and the component $K_{x}\subset\Omega$ is a $C^{1}-$ curve homeomorphic to the circle. We call such $K_{x}$ an admissible cycle.
\begin {lem}
Let $\Omega\subset\mathbb{R}^{2}$ be a bounded multiply connected domain with Lipschitz boundary $\partial\Omega=\cup^{N}_{j=1}\Gamma_{j}.$ Assume that  $\mathbf{v}, \mathbf{B}\in W^{1,2}(\Omega)$ and $p\in W^{1,s}(\Omega),s\in[1,2)$ satisfy the ideal MHD system (3.1) for almost all $ x\in\Omega.$ Assume that there exists a sequence of functions $\{\Phi_{\mu}\}$ such that $\Phi_{\mu}\in W_{loc}^{1,s}(\Omega)$ and $\Phi_{\mu}\rightharpoonup\Phi$ in $W_{loc}^{1,s}(\Omega)$ for some $s\in[\frac{4}{3},2).$ Then $\Phi_{\mu}|_{S}$ converges $\Phi|_{S}$ uniformly on almost all admissible cycles S. Moreover, if $\Phi_{\mu}$ satisfies the one-side maximum principle locally, then $\Phi$ satisfies the one-side maximum principle globally.
\end {lem}
\begin{rem} For any "self-closed" system, if its behavior is enough well, from the mathematical viewpoint, its regularity is enough high, then this system must have some form of the "maximum principle" on energy, and the above behavior (or regularity) must have a sharp low bound. And if worse (or lower than this value) then this "good" behavior, this system must have some  singularity or instability behavior.

 Indeed,
Lemma 3.7 states the fact that the ideal MHD system have some maximum principle, although it is a hyperbolic equations. from the analysis of the above paragraph,  it doesn't contradict with our known facts of the hyperbolic system.  In my opinion, The maximum principle is only a expression form of energy stability, then for any system in nature, if energy is a stability form in this system, we should find a form of  the maximum  principle of some variable in this system, of course including the hyperbolic system.
\end{rem}

$$$$

Although we do not know whether the function $\Phi$ is continuous or not on $\Omega,$ we shall prove that $\Phi$  has some continuity properties on stream lines.

Let $K\in T_{\psi}$ with $diam K>0.$ Take any $x\in K\setminus A_{w},$ and put $\Phi(K)=\Phi(x).$ This definition is correct by  Theorem 3.3.

\begin {lem}
Let $A,B\in T_{\psi},\textmd{diam} A>0,\textmd{diam }B>0.$ Consider the corresponding arc $[A,B]\subset T_{\psi}$ joining A to B. Then the restriction $\Phi|_{[A,B]}$ is a continuous function.
\end {lem}

We say that a set $\mathcal{Z}\subset T_{\psi}$ has $ T-$measure zero if $\mathfrak{h}^{1}(\{\psi(C):C\in\mathcal{Z}\})=0.$
In the following, we show that the function $\Phi|_{T_{\psi}}$ has some analogs of Luzin's N-property.

\begin {lem}
Let $A,B\in T_{\psi},\textmd{diam} A>0,\textmd{diam} B>0.$ If $\mathcal{Z}\subset[A,B]$ has $ T-$measure zero, then $\mathfrak{h}^{1}(\{\Phi(C):C\in\mathcal{Z}\})=0.$
\end {lem}

\begin {cor}
If $A,B\in T_{\psi},\textmd{diam }A>0,\textmd{diam} B>0,$ then
$\mathfrak{h}^{1}(\{\Phi(C):C\in[A,B]$
$and $
$ C$
$is$
$ not$
$ a$
$ regular$
$ cycle\})=0.$
\end {cor}
\begin{center}
\item\section{ Some properties of the MHD as $\nu,\eta\rightarrow0$}
\end{center}

  As mentioned in the above section, we assume that $\nu=\eta\rightarrow0,$  and rewrite the MHD equations in the domain $\Omega$ of type (1.1),
  \begin{equation}
\left\{\begin {array}{l}
-\nu\triangle \mathbf{Z}_{\nu}^{+}+\mathbf{Z}_{\nu}^{-}\cdot\nabla \mathbf{Z}_{\nu}^{+}+\nabla p_{\nu}=0,\\
-\nu\triangle \mathbf{Z}_{\nu}^{-}+\mathbf{Z}_{\nu}^{+}\cdot\nabla \mathbf{Z}_{\nu}^{-}+\nabla p_{\nu}=0,\\
\nabla\cdot \mathbf{Z}_{\nu}^{+} =0, \nabla\cdot \mathbf{Z}_{\nu}^{-}=0,\\

 \mathbf{Z}_{\nu}^{+}(x)|_{\partial\Omega}=\mathbf{a}_{\nu}(x)+\mathbf{b}_{\nu}(x), \mathbf{Z}_{\nu}^{-}(x)|_{\partial\Omega}=\mathbf{a}_{\nu}(x)-\mathbf{b}_{\nu}(x),
\end{array}\right.
\end{equation}

where $\mathbf{a}_{\nu}(x)=\mathbf{b}_{\nu}(x)=o(\nu).$

Here we assume the norms $\|\mathbf{Z}^{\pm}_{\nu}\|_{W^{1,2}(\Omega)},$ and $\|p_{\nu}\|_{W^{1,s}(\Omega)}$ are uniformly bounded for each $s\in[1,2)$ and $\nu.$ Moreover,  $\mathbf{Z}^{\pm}_{\nu}\in W_{loc}^{3,2}(\Omega), p_{\nu}\in W_{loc}^{2,2}(\Omega).$
  Then
  $\mathbf{Z}^{\pm}_{\nu}\rightharpoonup \mathbf{w}^{\pm}$
   in
     $W^{1,2}(\Omega),\mathbf{B}_{\nu}\rightharpoonup \mathbf{B}$
     in
     $ W^{1,2}(\Omega), p_{\nu}\rightharpoonup p $
      in
     $ W^{1,s}(\Omega),$ and $(\mathbf{w}^{\pm},p)$ satisfies the following equations:
\begin{equation}
\left\{\begin {array}{l}
\mathbf{w}^{+}\cdot\nabla \mathbf{w}^{-} +\nabla p=0,\\
\mathbf{w}^{-}\cdot\nabla \mathbf{w}^{+} +\nabla p=0,\\
\nabla\cdot\mathbf{ w}^{+}=\nabla\cdot\mathbf{ w}^{-}=0,\\
\mathbf{w}^{+}|_{\partial\Omega}=0,\mathbf{w}^{-}|_{\partial\Omega}=0.
\end{array}\right.
\end{equation}

According to Theorem 3.3 in the above section, we set $\Phi^{\pm}=p+\frac{1}{2}| \mathbf{w}^{\pm}|^{2}.$  From Lemma 3.7, in turn, first we consider the maximum (or extreme ) points of $\Phi^{pm},$ then try to the corresponding property of $\Phi^{\pm}_{\nu}$ in these points, where $\Phi^{\pm}_{\nu}=p_{\nu}+\frac{1}{2}|\mathbf{Z}^{\pm}_{\nu}|.$ Without generality of loss, we only study $\Phi^{-}$ and the corresponding function $\Phi^{-}_{\nu}.$ We divide $\Phi^{-}$ into the two cases:

 (a)  The maximum of $\Phi^{-}$ is attained on the boundary $\partial\Omega:$
 \begin {equation}
 \max_{j=0,\ldots,N}p_{j}=\textmd{ess}\sup_{x\in\Omega}\Phi(x).
 \end {equation}

 (b)  The maximum of $\Phi^{-} $ is not attained on the boundary $\partial\Omega:$
 \begin {equation}
 \max_{j=0,\ldots,N}p_{j}<ess\sup_{x\in\Omega}\Phi(x).
 \end {equation}

   In the case (a),  adding a constant to the pressure in (4.3), we can assume, without loss of generality, that
 $
 \max_{j=0,\ldots,N}p_{j}=\textmd{ess}\sup_{x\in\Omega}\Phi(x)=0.
$
In particular,
$
\Phi(x)\leq0$ in $\Omega.$

 Change(if necessary)  the numbering of the boundary components $\Gamma_{0},\Gamma_{1},\cdots,$
 $\Gamma_{N}$ in such a way that
  $
 p_{j}<0,\quad j=0,\ldots,M, \quad p_{M+1}=p_{M+2}\cdots=p_{N}=0,
$
 where $ M\in \{0,\ldots,N-1\}.$

 It is easy to compute directly that $\Phi^{-}_{\nu}$ satisfies the linear elliptic equation
 \begin {equation}
 \triangle\Phi^{-}_{\nu}=\frac{1}{\nu}\nabla\cdot(\Phi^{-}_{\nu}\mathbf{Z}^{-}_{\nu})+|\textmd{rot}\mathbf{Z}^{-}_{\nu}|^{2}
 -\nabla\cdot(2\mathbf{B}_{\nu}\cdot\nabla\mathbf{Z}^{-}_{\nu})
 +\frac{2}{\nu}\mathbf{B}_{\nu}\cdot\nabla\mathbf{Z}^{-}_{\nu}\cdot\mathbf{Z}^{-}_{\nu}.
 \end {equation}
 The key theorem  in this section is the following estimate.

\begin {thm} For any $i\in\mathbb{N},$ there exists $ \nu(i)\in\mathbb{N},$ such that the equality

$$\int_{S_{i\nu}(t)}|\nabla\Phi_{\nu}(x)|ds<C\bigg[\frac{1}{\ln \frac{1}{\nu}}+\frac{t}{\nu}\int_{S_{i\nu}(t)}\mathbf{Z}^{-}_{\nu}\cdot\mathbf{n}ds
$$
\begin {equation}
+2\bigg(\int_{S_{i\nu}(t)}|\mathbf{B}_{\nu}(x)|^{2}ds)^{\frac{1}{2}}(\int_{S_{i\nu}(t)}|\nabla\mathbf{Z}^{-}_{\nu}(x)|^{2}ds\bigg)^{\frac{1}{2}}
\end {equation}
$$+\frac{2}{\nu}\bigg(\int_{S_{i\nu}(t)}|\mathbf{B}_{\nu}(x)|^{2}ds\bigg)^{\frac{1}{2}}\bigg(\int_{S_{i\nu}(t)}|\mathbf{Z}^{-}_{\nu}(x)|^{4}ds\bigg)^{\frac{1}{2}}\bigg]$$
holds for every $\nu\geq \nu(i)$ and for almost all $t\in[\alpha_{2}t_{i},\alpha_{1}t_{i}],$ where the constant $ C $ is independent of $ t,i,\nu.$
\end {thm}
\begin{rem}   Because the ideal MHD system satisfies the Bernoulli law from Theorem 3.3, it can be regarded as the equilibrium.   The above result in Theorem 4.1 implies that the total energy nearby the boundary in MHD system in (4.1) is a stability state on the total energy near the corresponding ideal MHD system. At the same time, the above inequality implies that the MHD system in (4.1) preserves the one-side maximum  principle locally if   the corresponding ideal MHD  has one-side maximum principle globally in the domain.
\end{rem}

To prove Theorem 4.1, according to the analysis in Remark 4.2,  in fact, we have to prove that the MHD system (4.1) has a one-side maximum principle locally since (4.1) has a high nonlinear interaction such that it may not have maximum principle globally. Therefore , we have to construct the approximate subsets of the preimage of the corresponding local maximum value. Here we still decompose it into several steps.

$\mathbf{4.1}$ Construct the approximate integral line $S_{i\nu}.$

Denote by $B_{0},\ldots,B_{N},$ the elements of $T_{\psi}$ such that $B_{j}\supset\Gamma_{j},j=0,\ldots,N.$   Then every element $C\in[B_{i},B_{j}]\setminus\{B_{i},B_{j}\}$  is a connected component of a level set of $\psi$ such that the sets $B_{i}, B_{j}$ lie in different connected components of $\mathbb{R}^{2}\setminus C. $ Put $\alpha=\max_{j=0,\ldots,M}\min_{C\in[B_{i},B_{j}]}\Phi(C),$ then $\alpha<0.$

Step 1. Choose a sequent regular values $\{t_{i}\}^{\infty}_{i=1}$ satisfying the following properties:(a),$\{t_{i}\}^{\infty}_{i=1}\subset(0,-\alpha);$ (b),
$t_{i+1}=\frac{1}{2}t_{i};$ (c), if $\Phi(C_{i})=-t_{i},$ then $ C_{i} \in[B_{i},B_{N}]$ is a regular cycle.

Step 2. For $j=0,\ldots,M$ and $i\in\mathbb{N},$ put $A^{j}_{i}=\max\{C\in[B_{j},B_{N}]:\Phi(C)=-t_{i}\},$ then each set $A^{j}_{i}\subset\Omega$ is a $C^{1}-$curve homeomorphic to the unit circle. In particular, for each $i\in\mathbb{N},$ the compact set $\cup^{M}_{j=0}A^{j}_{i}$ is separated from $\partial\Omega$ and $\textmd{dist}(\cup^{M}_{j=0}A^{j}_{i},\partial\Omega)>0.$  Then for each $i$ and for sufficiently small $h(i)>0,$
 $\{x\in\Omega: \textmd{dist}(x,\Gamma_{h})<h\}\subset\Omega\setminus(\cup^{M}_{j=0}A^{j}_{i})$ is homeomorphic to the open ring.

 Step 3. For small $h,$ there exists some connected component $V_{i}$ such that $\{x\in\Omega: \textmd{dist}(x,\Gamma_{h})<h\}\subset V_{i}.$ At the same time, $V_{i}$ satisfies the following properties: (i),$(\Gamma_{N}\subset V_{i};$
(ii), $\Omega\cap\partial V_{i}=A^{0}_{i}\cup\ldots\cup A^{M}_{i};$ (iii),  the sequence $\{V_{i}\}$ is decreasing. From (iii), the sequence $\{\partial\Omega\cap(\partial V_{i})\}$ is nonincreasing,and for sufficiently large $i,$ the set $\partial\Omega\cap(\partial V_{i})$ is independent of $i.$
Set $\partial\Omega\cap(\partial V_{i})=\Gamma_{k}\cup\ldots\cup\Gamma_{N}.$ Therefore,  $\partial V_{i}=A^{0}_{i}\cup\ldots\cup A^{M}_{i}\cup\Gamma\ldots\cup\Gamma_{N}.$

Step 4. From Lemma 3.7, the uniform convergence $\Phi_{\nu}|_{A^{j}_{i}}\rightrightarrows\Phi|_{A^{j}_{i}}=-t_{i}$ as  $\nu\rightarrow0.$ Thus, for every $i\in\mathbb{N},$ there exists $\nu_{i}$ such that for all $\nu\leq\nu_{i},$
$\Phi_{\nu}|_{A^{j}_{i}}<-\frac{7}{8}t_{i},$ $\Phi_{\nu}|_{A^{j}_{i}}>-\frac{5}{8}t_{i},$ for any $j=0,\ldots,M.$

Step 5. Denote by $W^{j}_{i\nu}(t)$ the connected component of the open set $\{x\in V_{i}\setminus \bar{V}_{i+1}:\Phi_{\nu}(x)>-t\}$ such that $\partial W^{j}_{i\nu}(t)\supset A^{j}_{i+1},$ for $\nu\leq\nu(i),j=1,\ldots,M,$ and $t\in[\frac{5}{8}t_{i},\frac{7}{8}t_{i}].$

Step 6. Put $W_{i\nu}(t)=\cup^{M}_{j=0}W^{j}_{i\nu}(t),$
$ S_{i\nu}(t)=(\partial W_{i\nu}(t))\cap V_{i}\setminus \bar{V}_{i+1}.$ Then  the properties of $\Phi_{\nu}$ on $ S_{i\nu}$

 (i) $\Phi_{\nu}\equiv-t$ on $ S_{i\nu}(t).$ It is obvious from the construction.

 (ii)  $\Phi_{\nu}$ is differentiable (in classical sense) at every point $x\in S_{i\nu}(t)$ with $\nabla\Phi_{\nu}(x)\neq0.$
    And the values $ t\in[\alpha_{2}t_{i},\alpha_{1}t_{i}]$ having the above property will be called $(\nu,i)-regular. $

(iii) By construction, the following form holds,$
\int_{S_{i\nu}(t)}\nabla\Phi_{\nu}(x)\cdot \mathbf{n}ds=$
$-\int_{S_{i\nu}(t)}$
$|\nabla\Phi_{\nu}(x)|ds<0,
$
where $\mathbf{n}$ is the unit outward (with respect to $  W_{i\nu}(t)$) normal vector to $\partial W_{i\nu}(t).$

$\mathbf{4.2}.$ The proof of Theorem 4.1.

$ Proof$.  First,  we claim $
 \int_{\Gamma_{h}}\Phi_{\nu}^{2}ds <\sigma^{2},\quad \forall h\in(0,\delta_{\sigma}],\forall k\geq k',
 $
 for any sufficiently small $\sigma>0.$

 Indeed, we fix $i\in\mathbb{N},$ and assume $\nu\leq \nu(i)$ and choose the parameter $\delta_{\sigma}\in(0,\delta_{0}]$ small enough to satisfy the following conditions:
$
\int_{\Gamma_{h}}\Phi^{2}ds<\frac{1}{3}\sigma^{2} \quad \forall h\in(0,\delta_{\sigma}],
$
Since for any $q\in(1,2),$ the norm $\|\Phi_{\nu}\|_{W^{1,q}(\Omega)}$ are uniformly bounded, then the norms $\|\Phi_{\nu}\nabla\Phi_{\nu}\|_{L_{loc}^{q}(\Omega)}$ are uniformly bounded as well. In particular, for $q=\frac{6}{5},$ we have
 $|\int_{\Gamma_{h'}}\Phi_{\nu}^{2}ds-\int_{\Gamma_{h''}}\Phi_{\nu}^{2}ds|
 \leq2\int_{\Omega_{h''}\setminus\Omega_{h'}}|\Phi_{\nu}\nabla\Phi_{\nu}|.$

(4.11) holds from the following inequality (4.13),
\begin {equation}
-\frac{1}{3}\sigma^{2} <\int_{\Gamma_{h'}}\Phi_{\nu}^{2}ds-\int_{\Gamma_{h''}}\Phi_{\nu}^{2}ds<\frac{1}{3}\sigma^{2}\quad\forall h',h''\in(0,\delta_{\sigma}],\forall \nu\in\mathbb{N}.
\end {equation}
$$\leq2(\int_{\Omega_{h''}\setminus\Omega_{h'}}|\Phi_{\nu}\nabla\Phi_{\nu}|^{\frac{6}{5}}dx)^{\frac{5}{6}}
\textmd{ meas}(\Omega_{h''}\setminus\Omega_{h'})^{\frac{1}{6}}\rightarrow0,\quad as\quad h',h''\rightarrow0.$$

  By direct computation,
  $$\int_{S}\nabla\Phi_{\nu}\cdot\mathbf{n}ds=\int_{S}\textmd{rot}\mathbf{Z}^{-}_{\nu}\cdot \mathbf{Z}^{-}_{\nu}\cdot\mathbf{n}ds-2\int_{S}\mathbf{B}_{\nu}\cdot\nabla \mathbf{Z}^{-}_{\nu}\cdot\mathbf{n}ds.$$

For small $h>0,$ denote $\Gamma_{h}=\{x\in\Omega: \textmd{dist}(x,\Gamma_{K}\cup\cdots\cup\Gamma_{N})=h\},\Omega_{h}=\{x\in\Omega: \textmd{dist}(x,\Gamma_{K}\cup\cdots\cup\Gamma_{N})<h\}.$ By elementary results of analysis, there is a positive constant
$\delta_{0}<\frac{1}{2}\min\{|x-y|:x\in\Gamma_{j},y\in\Gamma_{m},j,m\in\{0,\ldots,N\},j\neq m\}$
such that for each $ h\leq\delta_{0},$ the set $\Gamma_{h}$ is a union of  N-K+1 $ C^{1}-$smooth curve homeomorphic to the circle, and
$
\mathfrak{h}^{1}(\Gamma_{h})\leq C_{0},\quad \forall h\in(0,h_{0}],
$
where $ C_{0}=3\mathfrak{h}^{1}(\Gamma_{K}\cup\cdots\cup\Gamma_{N})$ is independent of $ h.$

  Now, fix a sufficiently small $\varepsilon>0.$ The exact value of $\varepsilon$ will be specified below. Our next purpose is as follows. For a given sufficiently large $\nu\geq \nu',$ and a number $h_{\nu}\in(0,\delta_{\sigma})$ such that the estimates

\begin {equation}
\bigg|\int_{\Gamma_{h_{\nu}}}\textmd{rot}\mathbf{Z}^{-}_{\nu}\cdot \mathbf{Z}^{-}_{\nu}\cdot\mathbf{n}ds\bigg|<\varepsilon,\quad\int_{\Gamma_{h_{\nu}}}|\mathbf{Z}^{\pm}_{\nu}|^{2}ds<C_{\varepsilon}\nu^{2},\quad \int_{\Gamma_{h_{\nu}}}|\mathbf{B}_{\nu}|^{2}ds<C_{\varepsilon}\nu^{2},
\end {equation}
where the constant $C_{\varepsilon}$ is independent of $ \sigma. $ Take $\Gamma=\Gamma_{K}\cup\cdots\cup\Gamma_{N},$ and consider the function $g(h)=\int_{\Gamma_{h}}|\mathbf{Z}^{\pm}_{\nu}|^{2}ds.$ In particular,
$g(0)=\int_{\Gamma}|\mathbf{Z}^{\pm}_{\nu}|^{2}ds=\|\mathbf{a}\pm \mathbf{b}\|^{2}_{L^{2}(\Gamma)}.$
 Also denote $f(h)=\int_{\Gamma_{h}}|\nabla\mathbf{Z}^{\pm}_{\nu}||\mathbf{Z}^{\pm}_{\nu}|ds.$ By the classical formula of changing variables in the integral, there exists a $ C^{1}-$smooth function $J:\bar{\Omega}_{\delta_{0}}\rightarrow(0,\infty)$(not depending on $ \nu$) such that $J_{\Gamma}\equiv1$ and \footnote {Here $J(x)$ is the Jacobian of the following mapping:$\varphi:\bar{\Omega}_{\delta_{0}}\ni x\mapsto\varphi(x)=(\gamma(x),dist(x,\Gamma))\in\Gamma\times[0,\delta_{0}],$ where $\gamma(x)\in\Gamma$ is a metric projection of $ x$ onto $\Gamma:|x-\gamma(x)|=dist(x,\Gamma).$}
 $\bigg(\int_{\Gamma_{h}}J|\mathbf{Z}^{\pm}_{\nu}|^{2}ds\bigg)'_{h}\leq2\int_{\Gamma_{h}}J|\nabla\mathbf{Z}^{\pm}_{\nu}||\mathbf{Z}^{\pm}_{\nu}|ds.$
 Consequently, there are constants $C_{1},C_{2}>0$ (not depending on $ h$) such that for every $h_{0}\in(0,\delta_{0}],$ the following estimate holds:
 \begin{equation}
 \ln\bigg(\frac{C_{1}g(h_{0})}{\nu^{2}}\bigg)\leq C_{2}\int^{h_{0}}_{0}\frac{f(h)}{g(h)}dh.
\end{equation}
 Put
 $C_{\varepsilon}=\frac{1}{C_{1}}\exp(\frac{2C_{2}}{\varepsilon}).$
 Consider two possible cases:

 Case 1.$ g(h)\leq C_{\varepsilon}\nu^{2},\forall h\in[0,\delta_{0}].$ Then by the H$\ddot{o}$lder inequality we obtain
 \begin {equation}
 \frac{1}{\delta_{\sigma}}\int_{\Omega_{\delta_{\sigma}}}|\nabla\mathbf{Z}^{\pm}_{\nu}||\mathbf{Z}^{\pm}_{\nu}|dx\leq \frac{\sqrt{\delta_{\sigma}C_{\varepsilon}\nu^{2}}}{\delta_{\sigma}}\bigg(\int_{\Omega_{\delta_{\sigma}}}
 |\nabla\hat{\mathbf{Z}}^{\pm}_{k}|^{2}\bigg)^{\frac{1}{2}}
 <\sqrt{\frac{2C_{\varepsilon}\nu^{2}}{\delta_{\sigma}}}.
 \end {equation}
 Here we have used the estimate
$
 \int_{\Omega}|\nabla\mathbf{Z}^{\pm}_{\nu}|^{2}dx<2,
 $
 which is valid for sufficiently small $ \nu $( because of the convergence $ \int_{\Omega}|\nabla Z^{\pm}_{\nu}|^{2}dx\rightarrow1$ ). The above statement  implies that there exists $\bar{h}_{\nu}\in(0,\delta_{\sigma})$ such that
$
 \int_{\Gamma_{\bar{h}_{\nu}}}|\mathbf{Z}^{\pm}_{\nu}||\nabla\mathbf{Z}^{\pm}_{\nu}|ds<\sqrt{\frac{2C_{\varepsilon}\nu^{2}}{\delta_{\sigma}}}.
$
 Then, taking into account that $\nu\rightarrow0,$ while $\delta_{\sigma}$ is independent of $ \nu,$ we obtain the required estimates (4.8) for sufficiently small $\nu.$

 $ Case 2. \quad\sup_{h\in[0,\delta_{\sigma}]}g(h)> C_{\varepsilon}\nu^{2}.$ Take $h_{0}=\min\{h\in[0,\delta_{\sigma}]:g(h)=C_{\varepsilon}\nu^{2}\}.$ By choice of $C_{\varepsilon}$ and (4.9), then
 $
 \frac{2}{\varepsilon}\leq\int^{h_{0}}_{0}\frac{f(h)}{g(h)}dh.
 $

 We claim that there exists $\bar{h}_{\nu}\in(0,h_{0})$ satisfying (4.8). Suppose the contrary;then $f(h)\geq\varepsilon$ for all $h\in(0,h_{0}).$ By the H$\ddot{o}$lder inequality,
 $f^{2}(h)\leq g(h)\int_{\Gamma_{\bar{h}_{\nu}}}|\nabla\mathbf{Z}^{\pm}_{\nu}|^{2}ds.$
 Consequently,
 $\int_{\Gamma_{\bar{h}_{\nu}}}|\nabla\mathbf{Z}^{\pm}_{\nu}|^{2}ds\geq\frac{f^{2}(h)}{g(h)}\geq\frac{f(h)\varepsilon}{g(h)},\quad h\in(0,h_{0}). $
 Hence, $$\int_{\Omega_{h_{0}}}|\nabla\mathbf{Z}^{\pm}_{\nu}|^{2}dx=\int^{h_{0}}_{0}dh\int_{\Gamma_{\bar{h}_{\nu}}}
 |\nabla\mathbf{Z}^{\pm}_{\nu}|^{2}ds\geq\int^{h_{0}}_{0}\frac{f(h)\varepsilon}{g(h)}dh\geq2.$$
 ( In the last inequality we have used the well-known identity $ |\nabla \textmd{dist}(x,\Gamma)|\equiv 1$ on $\Omega_{h_{0}}.$)We have obtained the contradiction. This proves the existence of the required $\bar{h}_{\nu}\in(0,\delta_{\sigma})$ satisfying (4.8) for sufficiently small $ \nu.$

Now, for $ (\nu,i)-$regular value $t\in[\alpha_{2}t_{i},\alpha_{1}t_{i}],$ consider the domain
$\Omega_{i\bar{h}_{\nu}}(t)=W_{i\nu}(t)\cup\bar{V_{i+1}}\setminus\bar{\Omega}_{\bar{h}_{\nu}}.$
Since $
\partial W_{i\nu}(t)=S_{i\nu}(t)\cup A^{0}_{i+1}\cdots\cup A^{M}_{i+1},$
 $\partial\Omega_{i\bar{h}_{\nu}}(t)=\Gamma_{\bar{h}_{\nu}}\cup S_{i\nu}(t).$ Integrating the equation
$
 \triangle\Phi^{-}_{\nu}=\frac{1}{\nu}\nabla\cdot(\Phi^{-}_{\nu}\mathbf{Z}^{-}_{\nu})+|\textmd{rot}\mathbf{Z}^{-}_{\nu}|^{2}
 -\nabla\cdot(2\mathbf{B}_{\nu}\cdot\nabla\mathbf{Z}^{-}_{\nu})
 +\frac{2}{\nu}\mathbf{B}_{\nu}\cdot\nabla\mathbf{Z}^{-}_{\nu}\cdot\mathbf{Z}^{-}_{\nu}
 $
over the domain $\Omega_{i\bar{h}_{\nu}}(t),$ we have
$$\int_{S_{i\nu}(t)}\nabla\Phi_{\nu}\cdot\mathbf{n}ds+\int_{\Gamma_{\bar{h}_{\nu}}}\nabla\Phi_{\nu}\cdot\mathbf{n}ds$$
$$=\int_{\Omega_{\bar{h}_{\nu}}}|\textmd{rot}\mathbf{Z}^{-}_{\nu}|^{2}dx-2\int_{\Omega_{\bar{h}_{\nu}}}\textmd{div}(\mathbf{B}_{\nu}
\cdot\nabla\mathbf{Z}^{-}_{\nu})dx
+\frac{2}{\nu}\int_{\Omega_{\bar{h}_{\nu}}}\mathbf{B}_{\nu}\cdot\nabla\mathbf{Z}^{-}_{\nu}\cdot\mathbf{Z}^{-}_{\nu}dx$$
$$+\nu^{-1}\int_{S_{ik}(t)}\Phi_{\nu}\mathbf{Z}^{-}_{\nu}\cdot\mathbf{n}ds+\nu^{-1}\int_{\Gamma_{\bar{h}_{\nu}}}\Phi_{\nu}\hat{\mathbf{Z}}^{-}_{\nu}\cdot
\mathbf{n}ds$$
$$=\int_{\Omega_{\bar{h}_{\nu}}}|\textmd{rot}\mathbf{Z}^{-}_{\nu}|^{2}dx-2\int_{S_{i\nu}(t)}\mathbf{B}_{\nu}\cdot\nabla\mathbf{Z}^{-}_{\nu}\cdot\mathbf{n}ds
-2\int_{\Gamma_{\bar{h}_{\nu}}}\hat{\mathbf{B}}_{\nu}\cdot\nabla\mathbf{Z}^{-}_{\nu}\cdot\mathbf{n}ds$$
$$+2\nu^{-1}\int_{\Gamma_{\bar{h}_{\nu}}}\mathbf{B}_{\nu}|\mathbf{Z}^{-}_{\nu}|^{2}\cdot\mathbf{n}ds
+2\nu^{-1}\int_{S_{i\nu}(t)}\mathbf{B}_{\nu}|\mathbf{Z}^{-}_{\nu}|^{2}\cdot\mathbf{n}ds $$
$$+\nu^{-1}\int_{S_{i\nu}(t)}\Phi_{\nu}\mathbf{Z}^{-}_{\nu}\cdot\mathbf{n}ds+\nu^{-1}\int_{\Gamma_{\bar{h}_{\nu}}}\Phi_{\nu}\mathbf{Z}^{-}_{\nu}\cdot\mathbf{n}ds.$$

Furthermore,$$\int_{S_{i\nu}(t)}|\nabla\Phi_{\nu}|ds=\int_{\Gamma_{\bar{h}_{\nu}}}\nabla\Phi_{\nu}\cdot\mathbf{n}ds
-\int_{\Omega_{\bar{h}_{\nu}}}|\textmd{rot}\mathbf{Z}^{-}_{\nu}|^{2}dx$$
$$+2\int_{S_{i\nu}(t)}\mathbf{\mathbf{B}}_{\nu}\cdot\nabla\mathbf{Z}^{-}_{\nu}\cdot\mathbf{n}ds
+2\int_{\Gamma_{\bar{h}_{\nu}}}\mathbf{B}_{\nu}\cdot\nabla\mathbf{Z}^{-}_{\nu}\cdot\mathbf{n}$$
$$-2\nu^{-1}\int_{\Gamma_{\bar{h}_{\nu}}}\mathbf{B}_{\nu}|\mathbf{Z}^{-}_{\nu}|^{2}\cdot\mathbf{n}ds
-2\nu^{-1}\int_{S_{i\nu}(t)}\mathbf{B}_{\nu}|\mathbf{Z}^{-}_{\nu}|^{2}\cdot\mathbf{n}ds$$
$$-\nu^{-1}\int_{S_{i\nu}(t)}\Phi_{\nu}\mathbf{Z}^{-}_{\nu}\cdot\mathbf{n}ds-\nu^{-1}\int_{\Gamma_{\bar{h}_{\nu}}}\Phi_{\nu}
\mathbf{Z}^{-}_{\nu}\cdot\mathbf{n}ds$$
$$\leq\varepsilon-\int_{\Omega_{\bar{h}_{\nu}}}|\textmd{rot}\mathbf{Z}^{-}_{\nu}|^{2}dx+2\bigg(\int_{S_{i\nu}(t)}|\mathbf{B}_{\nu}|^{2}ds\bigg)^{\frac{1}{2}}
\bigg(\int_{S_{i\nu}(t)}|\nabla\hat{Z}^{-}_{\nu}|^{2}ds\bigg)^{\frac{1}{2}}+C_{\varepsilon}\nu^{-1}$$
$$+2\nu^{-1}\bigg(\int_{S_{i\nu}(t)}|\mathbf{B}_{\nu}|^{2}ds\bigg)^{\frac{1}{2}}(\int_{S_{i\nu}(t)}
|\mathbf{Z}^{-}_{\nu}|^{4}ds\bigg)^{\frac{1}{2}}
+C_{\varepsilon}\nu^{\frac{1}{2}}$$
$$+t\nu^{-1}\int_{S_{i\nu}(t)}\mathbf{Z}^{-}_{\nu}\cdot\mathbf{n}ds
+\bigg(\int_{\Gamma_{\bar{h}_{\nu}}}|\Phi_{\nu}|^{2}dx)^{\frac{1}{2}}(\int_{\Gamma_{\bar{h}_{\nu}}}
|\mathbf{Z}^{-}_{\nu}|^{2}dx\bigg)^{\frac{1}{2}}$$

$$\leq C\bigg[\varepsilon+\bigg(\int_{S_{i\nu}(t)}|B_{\nu}|^{2}ds\bigg)^{\frac{1}{2}}
\bigg(\int_{S_{i\nu}(t)}|\nabla\mathbf{Z}^{-}_{\nu}|^{2}ds\bigg)^{\frac{1}{2}}+C_{\varepsilon}\nu\bigg]$$
$$+2\nu^{-1}\bigg(\int_{S_{i\nu}(t)}|\mathbf{B}_{\nu}|^{2}ds\bigg)^{\frac{1}{2}}\bigg(\int_{S_{i\nu}(t)}|\mathbf{Z}^{-}_{\nu}|^{4}ds\bigg)^{\frac{1}{2}}
+C_{\varepsilon}\nu^{\frac{1}{2}}+\sigma C^{\frac{1}{2}}_{\varepsilon}+Ct\nu^{-1}\int_{S_{i\nu}(t)}\Phi_{\nu}\mathbf{Z}^{-}_{\nu}\cdot\mathbf{n}ds,$$
where the last second inequality obtains from H$\ddot{o}$lder inequality and Interpolation inequality.Choosing $\varepsilon=\frac{8}{\ln \nu^{-1}},\sigma=\frac{\nu^{2}}{C^{\frac{1}{2}}_{\varepsilon}},$ the proof of Theorem 4.1 is finished.

 Case 2. The maximum of $\Phi$ is not attained at $\partial\Omega.$  Adding a constant to the pressure, we assume, without loss of generality, that
\begin {equation}
\max_{j=0,\ldots,N}\hat{p}_{j}<0<\textmd{ess}\sup_{x\in\Omega}\Phi(x).
\end {equation}
(Here we do not exclude the case $\textmd{ess}\sup_{x\in\Omega}\Phi(x)=\infty.)$ Denote $\sigma=\max_{j=0,\ldots,N}$
$\hat{p}_{j}$
$<0.$

As in the previous subsection, we consider the behavior of $\Phi$ on the tree $T_{\psi}.$
\begin {lem}
There exists $F\in T_{\psi}$ such that $diam F>0,F\cap\partial\Omega=\emptyset,$ and $\Phi(F)>\sigma.$
\end {lem}

The proof is found in Lemma 3.10 of [9]. Fix $ F $ from above lemma and consider the behavior of $\Phi$ on the arcs
$[B_{j},F],j=0,\ldots,N.$(Recall that by $B_{j},$ we denote the elements of $T_{\psi}$ such that $\Gamma_{j}\subset B_{j}.$)
The rest part of this subsection is similar to that of Case 1.  with the following difference: $ F $ now plays the role that was played before by $B_{N},$ and the calculations become easier since $ F $ lie strictly inside $\Omega.$ Here we omit the detailed proof.
\begin {lem} For any $i\in\mathbb{N},$ there exists $ \nu(i)\in\mathbb{N},$ such that the equality
$$\int_{S_{i\nu}(t)}|\nabla\Phi_{\nu}(x)|ds<C\bigg(t+2\bigg(\int_{S_{i\nu}(t)}|\mathbf{B}_{\nu}(x)|^{2}ds\bigg)^{\frac{1}{2}}
\bigg(\int_{S_{i\nu}(t)}|\nabla\mathbf{Z}^{-}_{\nu}(x)|^{2}ds\bigg)^{\frac{1}{2}}$$

\begin {equation}
+2\nu^{-1}\bigg(\int_{S_{i\nu}(t)}|\mathbf{B}_{\nu}(x)|^{2}ds\bigg)^{\frac{1}{2}}\bigg(\int_{S_{i\nu}(t)}|\mathbf{Z}^{-}_{\nu}(x)|^{4}ds\bigg)^{\frac{1}{2}})
\end {equation}
holds for every $\nu\leq \nu(i)$ and for almost all $t\in[\alpha_{2}t_{i},\alpha_{1}t_{i}],$ where the constant $ C $ is independent of $ t,i,\nu.$
\end {lem}

$$$$

\begin{center}
\item\section{The existence theorem for any fixed $\nu=\eta\neq0$}
\end{center}

According to the results of the above two sections, in this section we study the existence of the MHD system  for any fix $\nu=\eta\neq0.$

 Now consider the MHD problem (1.2) (we use the form (1.5) in this section)in a bounded domain$\Omega$ with type (1.1). If the boundary $\mathbf{a},\mathbf{b}\in W^{\frac{3}{2},2}(\partial\Omega)$ satisfy the conditions (1.3),(1.4),
then there exist solenoidal extensions $\mathbf{A},\mathbf{B}\in W^{2,2}(\Omega)$ of $ \mathbf{a} ,\mathbf{b}$ (see[12,13]). We can find a weak solution $\mathbf{U}^{\mp}\in W^{2,2}(\Omega)$ of the Stokes problem such that $\mathbf{U}^{\mp}-\mathbf{A}\pm \mathbf{B}\in H(\Omega)$ and
$
\int_{\Omega}\nabla \mathbf{U}^{\mp}\cdot\nabla\eta dx=0\quad\quad\forall\eta\in H(\Omega).
$
Moreover,$
\|\mathbf{U}^{\pm}\|_{W^{2,2}(\Omega)}\leq c\|\mathbf{a}\pm \mathbf{b}\|_{W^{\frac{3}{2},2}(\partial\Omega)}.
$

By a weak solution of problem$(1.5),$ we understand a function $\mathbf{w}^{\pm}$ such that $\mathbf{Z}^{\mp}, \mathbf{w}^{\pm}=\mathbf{Z}^{\pm}-\mathbf{U}^{\pm}\in H(\Omega)$ satisfies
$$
\int_{\Omega}\nabla\mathbf{ w}^{\pm}\cdot\nabla\mathbf{\eta} dx-\int_{\Omega}((\mathbf{w}^{\pm}+\mathbf{U}^{\pm}-2\mathbf{B})\cdot\nabla)\mathbf{\eta}\cdot \mathbf{w}^{\pm}dx-\int_{\Omega}((\mathbf{w}^{\pm}-2\mathbf{B})\cdot\nabla)\mathbf{\eta}\cdot \mathbf{U}^{\pm}dx $$
\begin {equation}
=\int_{\Omega}(\mathbf{U}^{\pm}\cdot\nabla)\mathbf{\eta}\cdot \mathbf{U}^{\pm}dx
\end {equation}
for any $\mathbf{\eta}\in H(\Omega).$

 Without general of loss, we assume $\|\mathbf{Z}^{+}\mathbf{-U}^{+}\|_{H(\Omega)}\geq\|\mathbf{Z}^{-}-\mathbf{U}^{-}\|_{H(\Omega)}.$ In the following we only need to  consider the case $\mathbf{Z}^{+}.$
According to Riesz's Theorem, there exists elements $A\mathbf{w},\mathbf{F}$ in $ H(\Omega)$ such that
$$\int_{\Omega}((\mathbf{w}+\mathbf{U}-2\mathbf{B})\cdot\nabla)\mathbf{\mathbf{\eta}}\cdot \mathbf{w}dx-\int_{\Omega}((\mathbf{w}-2\mathbf{B})\cdot\nabla)\mathbf{\eta}\cdot \mathbf{U}dx+\int_{\Omega}(\mathbf{U}\cdot\nabla)\mathbf{\eta}\cdot\mathbf{U}dx$$
$$=[A\mathbf{w}+\mathbf{F},\mathbf{\eta}]£¬$$
where $\mathbf{w}=\mathbf{w}^{+},\mathbf{U}=\mathbf{U}^{+}.$

To prove the existence of a weak solution to the Magnetohydrodynamic problem(1.5), it is sufficient to show that all possible solutions of the integral identity
$$\int_{\Omega}\nabla\mathbf{ w}\cdot\nabla\mathbf{\eta} dx-\lambda\int_{\Omega}((\mathbf{w}+\mathbf{U}-2\mathbf{B})\cdot\nabla)\mathbf{\eta}\cdot \mathbf{w}dx-\lambda\int_{\Omega}((\mathbf{w}-2\mathbf{B})\cdot\nabla)\mathbf{\eta}\cdot \mathbf{U}dx$$
\begin {equation}
=\lambda\int_{\Omega}(\mathbf{U}\cdot\nabla)\mathbf{\eta}\cdot\mathbf{ U}dx,\quad \forall\eta\in H(\Omega)
\end {equation}
are uniformly bounded (with respect to $ \lambda\in[0,1]$) in $ H(\Omega).$

 Assume the above conclusion is false, then there exist sequences
${\lambda_{k}}\subset[0,1]$ and ${\mathbf{w}_{k}}\in H(\Omega), k\in\mathbb{N}$ such that

$$\int_{\Omega}\nabla \mathbf{w}_{k}\cdot\nabla\mathbf{\eta} dx-\lambda_{k}\int_{\Omega}((\mathbf{w}_{k}+\mathbf{U}-2\mathbf{B}_{k})\cdot\nabla)\mathbf{\eta}\cdot \mathbf{w}_{k}dx-\lambda_{k}\int_{\Omega}((\mathbf{w}_{k}-2\mathbf{B}_{k})\cdot\nabla)\mathbf{\eta}\cdot\mathbf{ U}dx$$
\begin {equation}
=\lambda_{k}\int_{\Omega}(\mathbf{U}\cdot\nabla)\mathbf{\eta}\cdot \mathbf{U}dx\quad \forall\mathbf{\eta}\in H(\Omega),
\end {equation}
where $\mathbf{u}_{k}+\mathbf{B}_{k}=\mathbf{w}_{k}+\mathbf{U},$ and
\begin {equation}
\lim_{k\rightarrow\infty}\lambda_{k}=\lambda_{0}\in[0,1],\quad\lim_{k\rightarrow\infty}J_{k}=\lim_{k\rightarrow\infty}\|\mathbf{w}_{k}\|_{ H(\Omega)}=\infty,
\end {equation}
where $\mathbf{w}=\mathbf{w}^{+},\|\mathbf{w}^{-}_{k}\|_{ H(\Omega)}\leq\|\mathbf{w}^{+}_{k}\|_{ H(\Omega)}$,since $\|\mathbf{Z}^{+}-\mathbf{U}^{+}\|_{H(\Omega)}\geq\|\mathbf{Z}^{-}-\mathbf{U}^{-}\|_{H(\Omega)}.$

It is easy to check that $\|\mathbf{B}_{k}\|_{ H(\Omega)}\leq C\|\mathbf{w}_{k}\|_{ H(\Omega)}$ where C is independent of $k.$

Let us take in (5.3) $\mathbf{\eta}=J^{-2}_{k}\mathbf{w}_{k}$ and denote $\widehat{\mathbf{w}}_{k}=J^{-1}_{k}\mathbf{w}_{k},\widehat{\mathbf{B}}_{k}=J^{-1}_{k}B_{k}.$ Since
$\int_{\Omega}((\mathbf{w}_{k}+\mathbf{U}-2\mathbf{B}_{k})\cdot\nabla)\mathbf{w}_{k}\cdot \mathbf{w}_{k}dx=0,$
we get
$$\int_{\Omega}|\nabla\widehat{ \mathbf{w}}_{k}|^{2}dx-\lambda_{k}\int_{\Omega}(\widehat{ \mathbf{w}}_{k}\cdot\nabla)\widehat{\mathbf{ w}}_{k}\cdot \mathbf{U }dx+2\lambda_{k}\int_{\Omega}(\widehat{\mathbf{B}}_{k}\cdot\nabla)\widehat{ \mathbf{w}}_{k}\cdot \mathbf{U} dx $$
\begin {equation}
=\lambda_{k}\int_{\Omega}(\mathbf{U}\cdot\nabla)\widehat{ \mathbf{w}}_{k}\cdot \mathbf{U}J^{-1}_{k}dx
\end {equation}
Since $\|\widehat{ \mathbf{w}}_{k}\|_{H(\Omega)}=1,\|\widehat{ \mathbf{B}}_{k}\|_{H(\Omega)}\leq C,$ there exist subsequences $\{\widehat{ \mathbf{w}}_{k_{l}}\},\{\widehat{ \mathbf{B}}_{k_{l}}\}$ converging weakly in $H(\Omega)$ to vector fields $\widehat{\mathbf{w}},\widehat{\mathbf{B}}\in H(\Omega).$  Because of the compact embedding
$$H(\Omega)\hookrightarrow L^{r}(\Omega)\quad\quad\forall r\in(1,\infty),$$
the subsequences $\{\widehat{ \mathbf{w}}_{k_{l}}\},\{\widehat{ \mathbf{B}}_{k_{l}}\}$ converge strongly in $ L^{r}(\Omega).$ Therefore, passing to a limit as $k_{l}\rightarrow\infty$ in(5.5), we obtain
\begin {equation}
1=\lambda_{0}\int_{\Omega}((\widehat{\mathbf{w}}-2\widehat{\mathbf{B}})\cdot\nabla)\widehat{\mathbf{w}}\cdot \mathbf{U} dx .
\end {equation}

Let us return to integral identity (5.3). Consider the functional
$$R_{k}(\mathbf{\eta})=\int_{\Omega}\nabla \mathbf{w}_{k}\cdot\nabla\mathbf{\mathbf{\eta}} dx-\lambda_{k}\int_{\Omega}(( \mathbf{w}_{k}+\mathbf{U}-2\mathbf{B}_{k})\cdot\nabla)\mathbf{\eta}\cdot\mathbf{ w}_{k} dx-\lambda_{k}\int_{\Omega}( \mathbf{w}_{k}\cdot\nabla)\mathbf{\eta}\cdot \mathbf{U }dx$$
$$+2\lambda_{k}\int_{\Omega}(\mathbf{B}_{k}\cdot\nabla)\mathbf{\eta}\cdot\mathbf{ U} -\lambda_{k}\int_{\Omega}(\mathbf{U}\cdot\nabla)\mathbf{\eta}\cdot\mathbf{ U}\quad\quad\mathbf{ \eta}\in\dot{H}^{1}(\Omega).$$
Obviously, $R_{k}(\mathbf{\eta})$ is a linear functional, and we have the estimate
$$|R_{k}(\mathbf{\eta})|\leq C(\|\mathbf{w}_{k}\|_{H(\Omega)}+\|\mathbf{w}_{k}\|^{2}_{H(\Omega)}+\|\mathbf{a}+\mathbf{b}\|^{2}_{W^{\frac{1}{2},2}(\partial\Omega)})\|\eta\|_{H(\Omega)},$$
with constant c independent of $k.$ It follows from (5.3) that
$$R_{k}(\mathbf{\eta})=0\quad\quad\forall\eta\in H(\Omega).$$
Therefore, there exist functions $p_{k}\in\hat{L}^{2}(\Omega)=\{q\in L^{2}(\Omega):\int_{\Omega}q(x)dx=0\}$ such that
$$R_{k}(\mathbf{\eta})=\int_{\Omega}p_{k}\nabla\cdot\mathbf{\eta} dx \quad\quad\forall\mathbf{\eta}\in \dot{H}^{1}(\Omega)$$
and \begin {equation}
\|p_{k}\|_{L^{2}(\Omega)}\leq C(\|\mathbf{w}_{k}\|_{H(\Omega)}+\|\mathbf{w}_{k}\|^{2}_{H(\Omega)}+\|\mathbf{a}+\mathbf{b}\|^{2}_{W^{\frac{1}{2},2}(\partial\Omega)}).
\end {equation}
The pair $(\mathbf{w}_{k},p_{k})$ satisfies the integral identity
$$\int_{\Omega}\nabla \mathbf{w}_{k}\cdot\nabla\mathbf{\eta} dx-\lambda_{k}\int_{\Omega}((\mathbf{w}_{k}+\mathbf{U}-2\mathbf{B}_{k})\cdot\nabla)\mathbf{\eta}\cdot \mathbf{w}_{k}dx-\lambda_{k}\int_{\Omega}((\mathbf{w}_{k}-2\mathbf{B}_{k})\cdot\nabla)\mathbf{\eta}\cdot\mathbf{ U}dx$$
\begin {equation}
-\lambda_{k}\int_{\Omega}(\mathbf{U}\cdot\nabla)\mathbf{\eta}\cdot \mathbf{U}dx=\int_{\Omega}p_{k}\nabla\cdot\mathbf{\eta} dx \quad \forall\mathbf{\eta}\in  \dot{H}^{1}(\Omega),
\end {equation}
Let $\mathbf{Z}^{+}_{k}=\mathbf{w}_{k}+\mathbf{U}.$ Then the identity (5.8) takes the form
$$\int_{\Omega}\nabla \mathbf{Z}^{+}_{k}\cdot\nabla\mathbf{\eta} dx-\int_{\Omega}p_{k}\nabla\cdot\mathbf{\eta} dx =-\lambda_{k}\int_{\Omega}((\mathbf{ Z}^{+}_{k}-2\mathbf{B}_{k})\cdot\nabla) \mathbf{Z}^{+}_{k}\cdot\eta dx\quad \forall\mathbf{\eta}\in  \dot{H}^{1}(\Omega).$$
Thus,$(\mathbf{Z}^{+}_{k},p_{k})$ might be considered as a weak solution to the Stokes problem in $\Omega$
\begin{equation}
\left\{\begin {array}{l}
-\triangle\mathbf{ Z}^{+}_{k}+\nabla p_{k}=\mathbf{f}_{k},\\

\nabla\cdot\mathbf{ Z}^{+}_{k} =0,\\
\mathbf{Z}^{+}_{k}|_{\partial\Omega}=\mathbf{a}(x)+\mathbf{b}(x),

\end{array}\right.
\end{equation}
with the right-hand side $\mathbf{f}_{k}=-\lambda_{k}((\mathbf{ Z}^{+}_{k}-2\mathbf{B}_{k})\cdot\nabla) \mathbf{Z}^{+}_{k}.$ Obviously, $\mathbf{f}_{k}\in L^{s}(\Omega)$ for $s\in(1,2)$ and
$$\|\mathbf{f}_{k}\|_{L^{s}(\Omega)}\leq C\|(( \mathbf{Z}^{+}_{k}-2\mathbf{B}_{k})\cdot\nabla) \mathbf{Z}^{+}_{k}\|_{L^{s}(\Omega)}\leq C\| \mathbf{Z}^{+}_{k}-2\mathbf{B}_{k}\|_{L^{\frac{2s}{2-s}}}\|\nabla \mathbf{Z}^{+}_{k}\|_{L^{2}(\Omega)}$$
$$\leq C(\|\mathbf{w}_{k}\|^{2}_{H(\Omega)}+\|\mathbf{a}+\mathbf{b}\|^{2}_{W^{\frac{1}{2},2}(\partial\Omega)})$$
where C is independent of $ k $. By well known local regularity results for the Stokes system (see[15]) we have the estimate
$$\| p_{k}\|_{W^{1,s}(\Omega')}\leq C(\|\mathbf{f}_{k}\|_{L^{s}(\Omega)}+\|\mathbf{Z}^{+}_{k}\|_{W^{1,2}(\Omega)}+\| p_{k}\|_{L^{2}(\Omega)})$$
\begin {equation}
\leq C(\|\mathbf{w}_{k}\|^{2}_{H(\Omega)}+\|\mathbf{a}+\mathbf{b}\|^{2}_{W^{\frac{1}{2},2}(\partial\Omega)}
+\|\mathbf{w}_{k}\|_{H(\Omega)}+\|\mathbf{a}+\mathbf{b}\|_{W^{\frac{1}{2},2}(\partial\Omega)})
\end {equation}
holds, where $\Omega'$ is an arbitrary domain with $\bar{\Omega'}\subset\Omega$ and the constant $ C $ depends on $dist(\Omega',\partial\Omega)$ but not on $k.$
Denote $\hat{p}_{k}=J^{-2}_{k}p_{k}.$ It follows from (5.7) and (5.9) that
$\|\hat{p}_{k}\|_{L^{2}(\Omega)}\leq Const,|\hat{p}_{k}\|_{W^{1,s}(\Omega')}\leq Const$
for any $\bar{\Omega'}\subset\Omega$ and $s\in(1,2).$ Repeating the above proceed, we have the estimates
 $\widehat{\mathbf{Z}}^{\pm}_{k},\widehat{\mathbf{B}}_{k}\in W^{3,2}_{loc}(\Omega),\quad \hat{p}_{k}\in W^{3,2}_{loc}(\Omega).$
 Hence, from the sequence ${\hat{p}_{k_{l}}}$ we can extract a subsequence,still denoted by ${\hat{p}_{k_{l}}},$ which converges weakly in $\hat{L}^{2}(\Omega)$ and $W^{1,s}_{loc}(\Omega)$
to some function $\hat{p}\in W^{1,s}_{loc}(\Omega)\cap\hat{L}^{2}(\Omega).$ Let $\varphi\in C^{\infty}_{0}(\Omega),$ taking
in (5.8) $\mathbf{\eta}=J^{-2}_{k}\mathbf{\varphi}$ and letting $k_{l}\rightarrow\infty$ yields
$$-\lambda_{0}\int_{\Omega}(\hat{\mathbf{w}}-2\widehat{\mathbf{B}})\cdot\nabla\varphi\cdot\widehat{\mathbf{w}}dx=\int_{\Omega}\hat{p}\nabla\cdot\mathbf{\varphi} dx\quad\forall\mathbf{\varphi}\in C^{\infty}_{0}(\Omega).$$
Integrating by parts in the last equality, we derive
\begin {equation}
\lambda_{0}\int_{\Omega}(\widehat{\mathbf{w}}-2\widehat{\mathbf{B}})\cdot\nabla\widehat{\mathbf{w}}\cdot\mathbf{\varphi} dx=-\int_{\Omega}\nabla\hat{p}\cdot\mathbf{\varphi} dx\quad\forall\mathbf{\varphi}\in C^{\infty}_{0}(\Omega).
\end {equation}
Hence,  $(\widehat{\mathbf{w}},\hat{p})$ satisfies, for almost all $x\in\Omega,$ the ideal MHD equations
\begin{equation}
\left\{\begin {array}{l}
\lambda_{0}(\widehat{\mathbf{w}}-2\hat{\mathbf{B}})\cdot\nabla\widehat{\mathbf{w}}+\nabla\hat{p}=0,\\
\nabla\cdot\widehat{\mathbf{w}}=\nabla\cdot\widehat{\mathbf{B}}=0,\\
\widehat{\mathbf{w}}|_{\partial\Omega}=\widehat{\mathbf{B}}|_{\partial\Omega}=0.
\end{array}\right.
\end{equation}
and $\widehat{\mathbf{w}}\in W^{1,2}(\Omega).\widehat{\mathbf{B}}\in W^{1,2}(\Omega),\hat{p}\in W^{1,s}(\Omega),$ for $s\in[1,2).$
By Lemma 3.6, $\hat{p}$ is a constant on each component of $\partial\Omega$ almost everywhere.

From the system (5.9), we put $\widehat{\mathbf{Z}}^{+}_{k}=\frac{1}{J_{k}}\mathbf{Z}^{+}_{k},\widehat{\mathbf{B}}_{k}=\frac{1}{J_{k}}\mathbf{B}_{k},\hat{p}_{k}
=\frac{1}{\lambda_{k}J^{2}_{k}}p_{k}.$ Multiplying
identities (5.9) by $\frac{1}{\lambda_{k}J^{2}_{k}},$ then the pair $(\widehat{\mathbf{Z}}^{+}_{k},\hat{p}_{k})$ satisfies the following system
\begin{equation}
\left\{\begin {array}{l}
-\frac{1}{\lambda_{k}J_{k}}\triangle \widehat{\mathbf{Z}}^{+}_{k}+\nabla \hat{p}_{k}=-(( \widehat{\mathbf{Z}}^{+}_{k}-2\widehat{\mathbf{B}}_{k})\cdot\nabla) \widehat{\mathbf{Z}}^{+}_{k},\\

\nabla\cdot \widehat{\mathbf{Z}}^{+}_{k} =\nabla\cdot \widehat{\mathbf{B}}^{+}_{k} =0,\\
\widehat{\mathbf{Z}}^{+}_{k}|_{\partial\Omega}=\widehat{\mathbf{a}}_{k}(x)+\widehat{\mathbf{b}}_{k}(x),

\end{array}\right.
\end{equation}
where $\widehat{\mathbf{a}}_{k}(x)=\frac{1}{J_{k}}\mathbf{a},\widehat{\mathbf{b}}_{k}(x)=\frac{1}{J_{k}}\mathbf{b}.$
\begin {lem}
 From the above statement of (5.9)-(5.13),we have that the norms $\|\widehat{\mathbf{Z}}^{+}_{k}\|_{W^{1,2}(\Omega)},\|\widehat{\mathbf{B}}^{+}_{k}\|_{W^{1,2}(\Omega)}$ and $\|\hat{p}^{+}_{k}\|_{W^{1,s}(\Omega)}$ are uniformly bounded for each $s\in[1,2)$

 Moreover,  $\widehat{\mathbf{Z}}^{+}_{k}\in W_{loc}^{3,2}(\Omega),\widehat{\mathbf{B}}^{+}_{k}\in W_{loc}^{3,2}(\Omega),\hat{p}^{+}_{k}\in W_{loc}^{2,2}(\Omega)$
  and
  $$\widehat{\mathbf{Z}}^{+}_{k}\rightharpoonup \widehat{\mathbf{w}}\quad in \quad W^{1,2}(\Omega),\widehat{\mathbf{B}}^{+}_{k}\rightharpoonup \widehat{\mathbf{B}}\quad in \quad W^{1,2}(\Omega),\hat{p}_{k}\rightharpoonup \hat{p}\quad in \quad W^{1,s}(\Omega).$$

 \end {lem}

Now, we receive the required contradiction using the Coarea formula and the result of Theorem in Section 4.

\begin {prop}
Assume that $\Omega\subset\mathbb{R}^{2}$ is a bounded domain of type (1.1) with $C^{2}-$smooth boundary $\partial\Omega,$ and $\mathbf{a},\mathbf{b}\in W^{\frac{3}{2},2}(\partial\Omega)$ satisfies the conditions (1.3),(1.4). Then Lemma 5.1 and (5.12) lead to a contradiction.
\end {prop}

Proof. For $i\in\mathbb{N}$ and $k\geq k(i)$ (see $\nu_{i}$ in section 4 ), put
$E_{i}=\cup_{t\in[\alpha_{2}t_{i},\alpha_{1}t_{i}]}S_{ik}(t).$
By the Coarea formula , for any integrable function $g:E_{i}\rightarrow\mathbb{R},$ the equality
\begin {equation}
\int_{E_{i}}g|\nabla\Phi_{k}|dx=\int^{\alpha_{1}t_{i}}_{\alpha_{2}t_{i}}\int_{S_{ik}(t)}g(x)d\mathfrak{h}^{1}(x)dt.
\end {equation}
holds. In particular, taking $g=|\nabla\Phi_{k}|$ and using(5.14), we obtain
$$\int_{E_{i}}|\nabla\Phi_{k}|^{2}dx=\int^{\alpha_{1}t_{i}}_{\alpha_{2}t_{i}}\int_{S_{ik}(t)}|\nabla\Phi_{k}|d\mathfrak{h}^{1}(x)dt\leq$$
$$\leq\int^{\alpha_{1}t_{i}}_{\alpha_{2}t_{i}}\bigg[\frac{1}{\ln J_{k}}+t+2\bigg(\int_{S_{ik}(t)}|\widehat{\mathbf{B}}_{k}(x)|^{2}ds\bigg)^{\frac{1}{2}}
\bigg(\int_{S_{ik}(t)}|\nabla\widehat{\mathbf{Z}}^{-}_{k}(x)|^{2}ds\bigg)^{\frac{1}{2}}$$
\begin {equation}
+2J_{k}\bigg(\int_{S_{ik}(t)}|\widehat{\mathbf{B}}_{k}(x)|^{2}ds\bigg)^{\frac{1}{2}}\bigg(\int_{S_{ik}(t)}|\widehat{\mathbf{Z}}^{-}_{k}(x)|^{4}ds
\bigg)^{\frac{1}{2}}\bigg]dt.
\end {equation}
$$\leq C\bigg(\frac{1}{\ln J_{k}}+t_{i}\bigg)+2\bigg(\int^{\alpha_{1}t_{i}}_{\alpha_{2}t_{i}}
\int_{S_{ik}(t)}|\widehat{\mathbf{B}}_{k}(x)|^{2}ds\bigg)^{\frac{1}{2}}\bigg(\int^{\alpha_{1}t_{i}}_{\alpha_{2}t_{i}}
\int_{S_{ik}(t)}|\nabla\widehat{\mathbf{Z}}^{-}_{k}(x)|^{2}ds\bigg)^{\frac{1}{2}}$$
$$+2J_{k}\bigg(\int^{\alpha_{1}t_{i}}_{\alpha_{2}t_{i}}\int_{S_{ik}(t)}|\widehat{\mathbf{B}}_{k}(x)|^{2}dsdt\bigg)^{\frac{1}{2}}
\bigg(\int^{\alpha_{1}t_{i}}_{\alpha_{2}t_{i}}\int_{S_{ik}(t)}|\widehat{\mathbf{Z}}^{-}_{k}(x)|^{4}dsdt\bigg)^{\frac{1}{2}}.$$
Now, taking $g=1$ in (5.14) and using the H$\ddot{o}$lder inequality, we have
\begin {equation}
\int^{\alpha_{1}t_{i}}_{\alpha_{2}t_{i}}\mathfrak{h}^{1}(S_{ik}(t))dt=
\int_{E_{i}}|\nabla\Phi_{k}|dx\leq\bigg(\int_{E_{i}}|\nabla\Phi_{k}|^{2}dx\bigg)^{\frac{1}{2}}
(\textmd{meas}(E_{i}))^{\frac{1}{2}}
\end {equation}
By construction, for almost all $t\in[\alpha_{2}t_{i},\alpha_{1}t_{i}],$ the set $S_{ik}(t)$ is the finite union of smooth cycles and $S_{ik}(t)$ separates $\Gamma_{j}$ from $\Gamma_{N}.$ In particular,
$$\mathfrak{h}^{1}(S_{ik}(t))\geq\min(\textmd{diam(}\Gamma_{j}),\textmd{diam}(\Gamma_{N})).$$

Hence, the left integral in (5.16) is greater than $Ct_{i},$ where $ C>0 $ does not depend on $ i .$ On the other hand, since
$t_{i}=O(\max(\|\mathbf{a}_{k}\pm \mathbf{b}_{k}\|^{\frac{1}{4}}_{W^{1,2}(\partial\Omega)},(\frac{1}{\ln J_{k}})^{\epsilon})),$ for sufficiently large $ i,k,$ the right integral in (5.16) is smaller than $Ct_{i}(\textmd{meas}(E_{i}))^{\frac{1}{2}}.$ However, $(\textmd{meas}(E_{i}))^{\frac{1}{2}}\rightarrow0$ as $i\rightarrow\infty.$ The obtained contradiction finishes the proof of Proposition 5.2.

$$$$

Proof of Theorem 1.2. Assume the conditions of Theorem 1.2 be satisfied. Suppose that its conclusion fails. Then, Lemma 5.1, there exist $\widehat{\mathbf{w}}^{\pm},\widehat{\mathbf{B}},\hat{p}$ and a sequence $(\widehat{\mathbf{Z}}^{\pm}_{k},\hat{p}_{k})$ satisfying the conditions of Lemma 5.1, and by Proposition 5.2, these assumptions lead to a contradiction.


$$$$

\end {document}